\documentstyle[twoside]{article}
\begin{document}
\newcommand{\ve}{\vspace{3mm}}
\headsep0.5cm
\textheight17.2cm
\textwidth11.5cm
\pagestyle{myheadings}
\evensidemargin1cm
\oddsidemargin1cm


\font\gordo=cmbx10 scaled \magstep2
\font\gordit=cmsl10 scaled \magstep1
\font\ggordo=cmbx10 scaled \magstep2
\font\ggordit=cmsl10 scaled \magstep1
\font\ninerm=cmr9
\font\tenrm=cmr10
\font\nineit=cmti9
\font\ninebf=cmbx9
\font\ninesl=cmsl9

 \newcommand{\be}{\begin{equation}}
\newcommand{\ee}{\end{equation}}
\newcommand{\bea}{\begin{eqnarray}}
\newcommand{\eea}{\end{eqnarray}}
\newcommand{\bean}{\begin{eqnarray*}}
\newcommand{\eean}{\end{eqnarray*}}

\newfont\frak{eufm10}

\ve
\parindent0cm

\newcommand{\R}{I\!\! R}

\newcommand{\mN}{I\!\!N}
\newcommand{\I}{1\!\mbox{I}}
\newcommand{\mQ}{\mbox{}\; l\!\!\!Q}
\newcommand{\Z}{Z\!\!\!Z}
\newcommand{\mC}{l\!\!\!C}
\newcommand{\mR}{I\!\!R}

\newtheorem{proposition}{Proposition}
\newtheorem{theorem}[proposition]{Theorem}

\newtheorem{corollary}[proposition]{Corollary}
\newtheorem{lemma}[proposition]{Lemma}

\newtheorem{definition}[proposition]{Definition}
\newtheorem{remark}[proposition]{Remark}

\newtheorem{example}[proposition]{Example}
\newtheorem{conclusion}[proposition]{Conclusion}

\newcounter{secnum}



\markboth{Bank, Giusti, Heintz, Mandel, Mbakop}{Polar Varieties}

\title{\bf Polar Varieties and Efficient Real Equation Solving: The Hypersurface
Case}

 \author{\sc B. Bank $^{1}$, M. Giusti $^{2}$, J. Heintz $^{3}$,\\[.5cm] 
                 \sc  R. Mandel, G. M. Mbakop $^{1}$}

 \maketitle

\addtocounter{footnote}{1}\footnotetext{Humboldt-Universit\"at zu Berlin,
Unter den Linden 6, D--10099 Berlin, bank@mathematik.hu-berlin.de; 
mbakop@mathematik.hu-berlin.de}

\addtocounter{footnote}{1}\footnotetext{GAGE, Centre de
Math\'ematiques. \'Ecole Polytechnique. F-91228, Palaiseau Cedex.
France,  giusti@ariana.polytechnique.fr}

\addtocounter{footnote}{1}\footnotetext{ Dept. de
Matem\'aticas, Estad\'{\i}stica y Computaci\'on. F. de Ciencias. U.
Cantabria. E-39071 Santander, Spain,  heintz@matsun1.unican.es}

\ve
\begin{abstract}
\noindent
The objective of this paper is to show how the recently
proposed method by Giusti, Heintz, Morais, Morgenstern, Pardo
\cite{gihemorpar} can be applied to a case of real polynomial equation
solving. Our main result concerns the problem of finding one
representative point for each connected component  of a real
bounded smooth hypersurface.\\ The algorithm in \cite{gihemorpar} yields a
method for symbolically solving a zero-dimensional polynomial
equation system in the affine (and toric) case. Its main feature
is the use of adapted data structure: Arithmetical networks and
straight-line programs.  The algorithm solves any affine
zero-dimensional equation system in non-uniform sequential time
that is polynomial in the length of the input description and an
adequately defined {\em affine degree} of the equation system.\\
Replacing the affine degree of the equation system by a suitably
defined {\em real degree} of certain polar varieties associated
to the input equation, which describes the hypersurface under
consideration, and using straight-line program codification of
the input and intermediate results, we obtain a method for the
problem introduced above that is polynomial in the input length
and the real degree.

\bigskip


\end{abstract}

\bigskip

\section{Introduction}

The present article is strongly related to the papers \cite{gihemorpar} and
\cite{gihemopar}.  Whereas the algorithms developed in these references are
related to the algebraically closed case, here we are concerned with the real
case.  Finding a real solution of a polynomial equation $f(x)=0$ where $f$ is a
polynomial of degree $d\ge 2$ with rational coefficients in $n$ variables is for
practical applications more important than the algebraically closed case.  Best
known complexity bounds for the problem we deal with are of the form $d^{O(n)}$
due to \cite{hroy}, \cite{rene}, \cite{basu}, \cite{sole}.  Related complexity
results can be found in \cite{canny}, \cite{grigo1}.  \par Solution methods for
the algebraically closed case are not applicable to real equation solving
normally.  The aim of this paper is to show that certain {\em polar varieties\/}
associated to an affine hypersurface possess a geometric invariant, {\em the
real degree\/}, which permits an adaptation of the algorithms designed in the
papers mentioned at the beginning.  The algorithms there are of "intrinsic
type", which means that they are able to distinguish between the semantical and
the syntactical character of the input system in order to profit both for the
improvement of the complexity estimates.  Both papers \cite{gihemorpar} and
\cite{gihemopar} show that the {\em affine degree\/} of an input system is
associated with the complexity when measured in terms of the number of
arithmetic operations.  Whereas the algorithms in \cite{gihemorpar} still need
algebraic parameters, those proposed in \cite{gihemopar} are completely
rational.  \par We will show that, under smoothness assumptions for the case of
finding a real zero of a polynomial equation of degree $d$ with rational
coefficients and $n$ variables, it is possible to design an algorithm of
intrinsic type using the same data structure, namely straight-line programs
without essential divisions and rational parameters for codifying the input
system, intermediate results and the output, and replacing the affine degree by
the real degree of the associated polar varieties to the input equation.  \par
The computation model we use will be an arithmetical network (compare to
\cite{gihemorpar}).  Our main result then consists in the following.  {\em There
is an arithmetical network of size $(nd\delta^*L)^{O(1)}$ with parameters in the
field of rational numbers which finds a representative real point in every
connected component of an affine variety given by a non-constant square-free
$n$-variate polynomial $f$ with rational coefficients and degree $d\ge 2$
$($supposing that the affine variety is smooth in all real points that are 
contained in it.$)$.  $L$
denotes the size of the straight-line program codifying the input and $\delta^*$
is the real degree associated to $f$.} \\
Close complexity results are the ones
following the approach initiated in \cite{ShSm93a}, and further developed in
\cite{ShSm93b}, \cite{ShSm93c}, \cite{ShSm93d}, \cite{ShSm1}, see also
\cite{Dedieu1}, \cite{Dedieu2}.  \par For more details we refer the reader to
\cite{gihemorpar} and \cite{gihemopar} and the references cited there.

\newpage

\section{Polar Varieties and Algorithms }
            

As usually, let $\mQ, \; \mR$ and $\mC$ denote the field of rational, real and 
complex numbers, respectively. The affine n--spaces over these fields are denoted by
$\mQ^n, \; \mR^n$ and $\mC^n$, respectively. Further, let $\mC^n$ be endowed
with the Zariski--topology, where a closed set consists of all common zeros of 
a finite number of polynomials with coefficients in $\mQ$.
Let $W \subset {\mC}^n$ be  a closed subset
with respect to this topology and let $W= C_1\cup\cdots \cup C_s$ be its 
decomposition into irreducible components with respect to the same topology. 
Thus $W, \; C_1,\ldots,C_s$ are algebraic subsets of  ${\mC}^n$. 
Let $1\le j \le s$, be arbitrarily fixed and consider the irreducible component
$C_j$ of $W$. 

In the following we need the notion of degree of an  affine algebraic variety.
Let $W \subset \mC^n$ be an algebraic subset given by a regular sequence 
$f_1, \cdots f_i \in \mQ[x_1, \cdots, x_n]$ of degree at most $d$. If 
$W \subset \mC^n$ ist zero--dimensional the {\it degree} of $W, \; deg W$, is
defined to be the number of points in $W$ (neither multiplicities nor points at
infinity are counted). If $W \subset \mC^n$ is of dimension greater than zero 
(i.e.
$dim W = n-i \ge 1$), then we consider the collection ${\cal M}$ of all affine varieties
of dimension $i$ given as the solution set in $\mC^n$ of a linear equation 
system $L_1 = 0, \; \cdots, L_{n-i} = 0$ with $L_{k} = \sum_{j=1}^{n} a_{kj} x_j
+ a_{k0}, \; a_{ki} \in \mQ, 1 \le i \le n$. Let ${\cal M}_{W}$ be the subcollection of
${\cal M}$ formed by all varieties $H \in {\cal M}$ such that the affine variety $H \cap W$
satisfies $H \cap W \not= \emptyset$ and $dim(H \cap W) = 0$. Then the affine 
degree of $W$ is defined as $max\{ \delta | \delta = deg(H \cap W), \;
H \in M_W \}$.

\begin{definition}\label{def1}
The component $C_j$ is called a {\rm real component} of $W$ if the real variety 
$C_j\cap \mR^n$ contains a smooth point of $C_j$. \par
\noindent If we denote 
\[
I = \{ j \in \mN | 1 \le j \le s,  \hskip 3pt C_j \hskip
3pt{\hbox {\rm  is a real component of $W$}} \}.
\] 
then the affine variety $W^\ast := \bigcup \limits_{j \in I} C_j \; \subset 
\mC^n$ 
{\it is called the real part} of $W$. 
By $deg^{\ast} W := degW^{\ast} = \sum\limits_{j \in I}deg C_j$ we define
the {\it real degree} of the set  $W$.
\end{definition}

\begin{remark}
{\rm Observe that $deg^{\ast} W= 0$ holds if and only if the real part $W^{\ast}$ of $W$ 
is empty.}
\end{remark}

\begin{proposition}\label{prop3}
  Let $f \in \mQ[X_1,\cdots,X_n]$ be a non-constant 
and square-free polynomial and let $\widetilde{V}(f)$
be the set of real zeros of the equation $f(x) =0$.
Assume $\widetilde{V}(f)$ to be bounded.
Furthermore,
let for every fixed $i, \; 0 \le i <n$, the real variety
 $$\widetilde{V}_i  :=  \{ x \in \mR^n | \; f(x) = {{\partial f(x)} \over {\partial X_1}} = 
 \ldots,{{\partial f(x)} \over {\partial X_i}} = 0 \}$$  
  be non-empty (and $\widetilde{V}_0$ is understood to be $\widetilde{V}(f)$). 
Suppose the variables to be in generic position. Then any point of $\widetilde{V}_i$
that is a smooth point of $\widetilde{V}(f)$ is also a smooth point of 
$\widetilde{V}_i$. Moreover, for every such point the Jacobian of the equation 
system $f=\frac{\partial f}{\partial X_1} = \cdots = \frac{\partial f}{\partial
X_i} =0$ has maximal rank.
\end{proposition}
\bigskip

{\bf Proof}

Consider the linear transformation $x \longleftarrow A^{(i)} y$, where 
the new variables are $y = (Y_1, \cdots, Y_n)$. Suppose that  
$A^{(i)}$ is given in the form
\[
\left( \begin{array}{ll}
I_{i,i} & 0_{i,n-i} \\
(a_{kl})_{n-i,i} & I_{n-i,n-i} \end{array} \right) ,
\]
where $I$ and 0 define a unit  and a zero matrix, respectively, and\\
$a_{kl} \in \R $ arbitrary if $k,l$ satisfy $i+1 \le k \le n, \;\;
1 \le l\le i$.\\
The transformation $x \longleftarrow A^{(i)} y$ defines a linear change of coordinates, since the
square matrix $A^{(i)}$ has full rank.

In the new coordinates, the variety $\widetilde{V}_i$ takes the form 
$$\widetilde{V}_i := \{ y \in \mR^n | \; f(y) = {{\partial f(y)} \over {\partial Y_1}} +
\sum_{j = i+1}^n a_{j1} {{\partial f(y)} \over {\partial Y_j}} =  \ldots =
{{\partial f(y)} \over {\partial Y_i}}  
+\sum_{j = i+1}^n a_{ji} {{\partial f(y)} \over {\partial Y_j}}= 0 \}$$

This transformation defines a map 
$\Phi_i \; : \mR^n \times \mR^{(n-i) i} \longrightarrow \mR^{i+1}$ given by
$$\Phi_i \left ( Y_1, \cdots, Y_i, \cdots, Y_n, a_{i+1, 1}, \cdots 
a_{n 1}, \cdots a_{i+1,i}, \cdots, a_{n, i} \right ) = $$
$$\left ( f,\; {{\partial f} \over {\partial Y_1}} +
\sum_{j = i+1}^n a_{j 1} {{\partial f} \over {\partial Y_j}},  \ldots,\;
{{\partial f} \over {\partial Y_i}}  
+\sum_{j = i+1}^n a_{j i} {{\partial f} \over {\partial Y_j}} \right )$$
For the moment let  
$$\alpha := (\alpha_1, \cdots, \alpha_{(n-i) i} ) :=
(Y_1, \cdots, Y_n, a_{i+1, \; 1}, \cdots a_{n,i}) \in \mR^n \times \mR^{(n-i) i}$$
Then the Jacobian matrix of $\Phi_i ( \alpha )$ is given by\newline
$J \left (\Phi_i (\alpha ) \right ) = \left ( 
{{\partial \Phi_i ( \alpha )} \over {\partial \alpha_j}} 
\right )_{(i+1)\times (n + (n-i) i) } = $
$$\left ( \begin{array}{ccllclccl}
 {{\partial f}\over{\partial Y_1}}& \cdots &
 {{\partial f} \over {\partial Y_n}}  & 0 & \cdots & 0 & \cdots & \cdots & 0 \\
 \ast & \cdots&   \ast  & {{\partial f} \over {\partial Y_{i+1}}} & \cdots & {{\partial f} \over {\partial Y_n}}
 & 0 \cdots  & \vdots & 0\\                 
 \vdots & & \vdots & \ddots & \ddots & 0 & \cdots & \ddots & 0 \\ 
 \ast & \cdots& \ast & 0 \cdots  & 0 \cdots & \cdots &{{\partial f} \over {\partial Y_{i+1}}}  
  & \cdots & {{\partial f} \over {\partial Y_n}}                  
                       \end{array}    \right ) $$
If $\alpha^0 = (Y_1^0, \cdots, Y_n^0, a_{i+1, \; 1}^0, \cdots a_{n, \; i}^0)$
belongs to the fibre $\Phi_i^{-1} (0)$, where $(Y_1^0, \cdots, Y_n^0)$ is a
point of the hypersurface $\widetilde{V} (f)$ and if there is an index 
$j \in \{ i+1, \cdots, n \}$ such that 
${{\partial f} \over {\partial Y_j}} \not= 0$ at this point, then the Jacobian matrix 
$J \left (\Phi_i (\alpha^0) \right )$ has the maximal rank $i + 1$.\\
Suppose now  that for all points of $\widetilde{V} (f)$ 
$${{\partial f(y)} \over {\partial Y_{i+1}}} = \cdots = 
{{\partial f(y)} \over {\partial Y_n}} = 0$$
and let $C := \mR^n \setminus \{ {{\partial f(y)} \over {\partial Y_1}} 
= \cdots = 
{{\partial f(y)} \over {\partial Y_n}} = 0 \}$, which is an open set. 
Then the restricted map 
\[
\Phi_i : C \times \R^{(n-i)i} \longrightarrow \R^{i+1}
\]
is transversal to the subvariety $\{ 0\}$ in $\R^{i+1}$.\\
By weak transversality due to {\sl Thom/Sard} (see e.g. \cite{golub}) 
applied to
the diagram

 $$ \begin{array}{lcc}
                    \Phi_i^{-1} (0)  & \hookrightarrow & \mR^n \times \mR^{(n-i) i} \\
                   & \searrow  & \downarrow \\
                      &  & \mR^{(n-i) i}
     \end{array}  $$

\noindent one concludes that the set of all $A \in \mR^{(n-i) i}$ for which transversality
holds is dense in $\mR^{(n-i) i}$.

Since the hypersurface $\widetilde{V}(f)$ is bounded by assumption, there is an open and
dense set of matrices $A$ such that the corresponding coordinate transformation 
leads to the desired smoothness. \hfill $\Box$\\

Let $f \in \mQ[X_1,\cdots,X_n]$ be a non--constant squarefree polynomial and let
$W := \{ x \in \mC^n | \; f(x) = 0 \}$ be the hypersurface defined by $f$. 
Consider the real variety $V := W \cap \mR^n$ and suppose:
\begin{itemize}
\item $V$ is non-empty and bounded,
\item the gradient of $f$ is different from zero in all points of $V$\\
(i.e. $V$ is a compact smooth hypersurface in $\mR^n$ and $f = 0$ is its 
regular equation)
\item the variables are in generic position.
\end{itemize}

\begin{definition}[Polar variety corresponding to a linear space]

Let $i, \; 0\leq i<n$, be arbitrarily fixed. Further, let 
$X^i := \{ x \in \mC^n | \; X_{i+1} = \cdots = X_n = 0 \}$ be the corresponding
linear subspace of $\mC^n$. Then, $W_i$ defined to be the Zariski closure of 
$$ \{ x \in \mC^n | \; f(x) = {{\partial f(x)} \over{\partial X_1}} = 
\cdots = {{\partial f(x)} \over {\partial X_i}} = 0, \; \Delta (x) :=
\sum_{j = 1}^{n} \left ( {{\partial f(x)} \over {\partial X_j}} \right )^2 
\not= 0 \} $$ 
is called the {\it polar variety} of $W$ associated to the linear 
subspace $X^i$. 
 
The corresponding real variety of $W_i$ is denoted by $V_i
:=  W_i \cap \mR^n$.
\end{definition}

\bigskip

\begin{remark}
{\rm
Because of the hypotheses that $V \not= \emptyset$ is a smooth hypersurface and
that $W_i \not= \emptyset$ by the assumptions above, the real variety $V_i
:=  W_i \cap \mR^n, \;
0 \le i <n$, is not empty and by smoothness of $V$, it has the description 
 $$V_i = \{ x \in \mR^n |  f(x) = {{\partial f(x)} \over {\partial X_1}} = 
 \ldots  = {{\partial f(x)} \over {\partial X_i}} = 0 \} .$$
($V_0$ is understood to be $V$.)\\
According to Proposition 3, $V_i$ is smooth if the coordinates are chosen to
be in generic position. Definition 4 of a polar variety is slightly different from the
one introduced  by L\^{e}/Teissier \cite{le}.  }
\end{remark}

\bigskip

\begin{theorem}

Let $f \in \mQ[X_1,\cdots,X_n]$ be a non--constant squarefree polynomial and let
$W :=  \{ x \in \mC^n | \; f(x) = 0 \}$ be the corresponding hypersurface.
Further, let $V := W \cap \mR^n$ be a non--empty, smooth, and bounded 
hypersurface in $\mR^n$ whose regular equation is given by $f = 0$. Assume the 
variables $X_1, \cdots, X_n$ to be generic. Finally, for 
every $i, \; 0 \le i < n$, let the polar varieties $W_i$ of $W$ corresponding
to the subspace $X^i$ be defined as above. Then it holds~:
\begin{itemize}
\item $V \subset W_0$, with $W_0 = W$ if and only if $f$ and $\Delta :=
\sum_{j=1}^n \left( \frac{\partial f}{\partial X_j} \right)^2$ are coprime,
\item $W_i$ is a non--empty equidimensional affine variety of dimension 
$n-(i+1)$ that is smooth in all its points that are smooth points of $W$,
\item the real part $W_i^\ast $ of the polar variety $W_i$ coincides with 
the Zariski closure in $\mC^n$ of 
$$V_i = \left\{ x \in \mR^n |  f(x) = 
{{\partial f(x)} \over {\partial X_1}} = \ldots = 
{{\partial f(x)} \over {\partial X_i}} = 0 \right\} ,$$
\item for any $j$, $i<j \le n$ the ideal  

      $$\left(f, {{\partial f} \over {\partial X_1}}, \ldots,{{\partial f}
      \over {\partial X_i}}\right)_{{{\partial f} \over {\partial X_j}}}$$ is
       radical.
\end{itemize}
\end{theorem}

{\bf Proof:}

Let $i, 0\le i < n$, be arbitrarily fixed. The first item is obvious
since $W_0$ is the union of all irreducible components of $W$ on which 
$\Delta$ does not vanish identically. \\
Then
$W_i$ is non-empty by the assumptions. The sequence $f, 
\frac{\partial f}{\partial X_1} , \ldots , \frac{\partial f}{\partial X_i}$ of 
polynomials of $\mQ[X_1,\ldots , X_n]$ forms a local regular sequence 
with respect to the smooth points of $W_i$ since
the affine varieties $\big\{ x \in \mC^n | f(x) = \frac{\partial f(x)}{\partial X_1}
= \cdots = \frac{\partial f(x)}{\partial X_k} = 0 \big\}$ and 
 $\big\{ x \in \mC^n | \frac{\partial f(x)}{\partial X_{k+1}} = 0\big\} $ are transversal
for any $k, 0\le k\le i-1$, by the generic choice of the coordinates, and hence
the sequence $f, \frac{\partial f}{\partial X_1} ,\cdots , \frac{\partial f}{\partial X_i}$
yields a local complete intersection with respect to the same points. 
This implies that $W_i$ are equidimensional 
and $dim_{\mC} W_i = n-(i+1)$ holds. We observe that every  smooth point of $W_i$
is a smooth point of $W$, which completes the proof of the second item.\\
The Zariski closure of $V_i$ is contained in $W_i^\ast$, which is a simple 
consequence of the smoothness of $V_i$. One obtains the reverse inclusion 
as follows. Let $x^\ast \in W_i^\ast$ be an arbitrary point, and let $C_{j\ast}$
be an irreducible component of $W_i^\ast$ containing this point, and $C_{j\ast}
\cap V_i \not= \emptyset$. Then
\[
\begin{array}{rl}
n-i-1 &= dim_{\R}(C_{j\ast}\cap V_i)= dim_{\R}
R(C_{j\ast}\cap V_i)=\\
&=dim_{\mC} R((C_{j\ast}\cap V_i)')\le dim_{\mC} C_{j\ast} = n-i-1,
\end{array}
\]
where $R(\cdot)$ and $(\,\, )'$ denote the corresponding sets of smooth points 
contained in $( \cdot )$ and the associated complexification, respectively. Therefore,
$dim_{\mC} (C_{j\ast} \cap V_i)' = dim_{\mC} C_{j\ast} = n-i-1$ and, hence,
 $C_{j\ast} =
(C_{j\ast}\cap V_i)'$, and the latter set is contained in the Zariski closure of $V_i$.\\
We  define the non-empty affine algebraic set 
\[
\widetilde{W}_i := \left\{ x\in \mC^n | f(x) =
\frac{\partial f(x)}{\partial X_1} = \cdots = \frac{\partial f(x)}{\partial 
X_i}
= 0 \right\} . 
\]
 Let $j, i<j\le n,$ be arbitrarily fixed. Then one finds a smooth 
point $x^\ast$ in $\widetilde{W}_i$ such that  $ \frac{\partial f(x^\ast)}{\partial X_j}
\not= 0$; let $x^\ast$ be fixed in that way. The hypersurface $W = \{ x \in \mC^n |
f(x) =0 \}$ contains $x^\ast$  as a smooth point, too.  Consider the local ring
${\cal O}_{W,x^\ast} $ of $x^\ast$  on the hypersurface $W$. (This is the ring
of germs of functions on $W$ that are regular at $x^\ast$. The local ring
${\cal O}_{W,x^\ast}$ is obtained by dividing the ring $\mbox{} \; \mC [X_1, \ldots , X_n ]$ of
polynomials by the principal ideal $(f)$, which defines $W$ as an affine variety, 
and then by localizing  at the maximal ideal $(X_1 - X_1^\ast,\ldots , X_n-X_n^\ast)$,
of the point $x^\ast = (X_1^\ast, \ldots , X_n^\ast)$ considered as a single
point affine variety.) Using now arguments from Commutative Algebra and Algebraic Geometry,
see e.g. Brodmann \cite{brod}, one arrives at the fact that ${\cal O}_{W,x^\ast}$ is an integral
regular local ring.\\
The integrality of  ${\cal O}_{W,x^\ast}$ implies that there is a uniquely determined irreducible
component $Y$ of $W$ containing the smooth point $x^\ast$ and locally this component
corresponds to the zero ideal of ${\cal O}_{W,x^\ast}$, which is radical. Since
the two varieties $\widetilde{W}_i \cap Y$ and $W\cap Y$ coincide locally, the 
variety $\widetilde{W}_i \cap Y$ corresponds locally to the same ideal.
 Thus, the desired radicality
is shown. This completes the proof.  \linebreak
$\mbox{} \hfill \Box$\\

\begin{remark}
{\rm If one localizes with respect to the function $ \Delta(x) = \sum\limits^n_{j=1}
\big( \frac{\partial f(x)}{\partial X_j}\big)^2$, then one obtains, in the same way
as shown in the proof above, that the ideal 
\[
\big( f, \frac{\partial f}{\partial X_1} , \ldots , \frac{\partial f}{\partial X_i}
\big)_\Delta
\]
is also radical.}
\end{remark}

\begin{remark}\label{rem8}
{\rm Under the assumptions of Theorem 6, for any $i,\;\; 0\le i<n$, we observe the 
following relations between the different non-empty varieties introduced up to now.
\[
V_i \subset V,\quad V_i \subset W^\ast_i \subset W_i \subset \widetilde{W}_i ,
\]
where $V$ is the considered real hypersurface, $V_i$ defined as in Remark 5,
$W_i$ the polar variety due to Definition 4, $W^\ast_i$ its real part according
to Definition 1, and $\widetilde{W}_i$ the affine variety introduced in the proof of 
Theorem 6. With respect to Theorem 6 our settings and assumptions imply that
$n-i-1 = dim_{\mC} \widetilde{W}_i = dim_{\mC} W_i = dim_{\mC} W^\ast_i =
dim_{\R} V_i $ holds. By our smoothness assumption and the generic choice of the 
variables we have for the respective sets of smooth points (denoted as before
by $R(\cdot ))$
\[
V_i = R(V_i)\subset R(W_i) \subset R(\widetilde{W}_i) \subset R(W) ,
\]
where $W$ is the affine hypersurface.\\

For the following we use the notations as before, fix an $i$ arbitrarily,$ \;\; 0\le i<n$, 
denote by $\delta^\ast_i$ the real degree of the polar variety $W_i$
(compare with Definition 1, by smoothness one has that the real degree of the 
polar variety  $W_i$ is equal to the real degree of the affine variety $\widetilde{W}_i$),
put $\delta^\ast := \max \{ \delta^\ast_k | 0 \le k \le i \}$ and let $ d := 
\deg f$. Finally, we write for shortness $r := n-i-1$.\\

We say that the variables $X_1,\ldots , X_n$ are in Noether position with respect
to a variety $\{ f_1= \cdots = f_s =0\}$ in $\mC^n, \; f_1, \ldots  , f_s \in 
\mQ [ X_1, \ldots , X_n],$ if, for each $r<k \le n$, there exists a polynomial of
$\mQ [X_1, \ldots , X_r, X_k]$ that is monic in $X_k$ and vanishes on 
$\{ f_1=\cdots =f_s=0\}.$ \\

Then one can state the next, technical lemma according to \cite{gihemorpar},
\cite{gihemopar},
 where the second reference is important in order to ensure that the occurring
straight-line programs use parameters in $\mQ$ only.}
\end{remark}
  
\begin{lemma}\label{lem9}
Let the assumptions of Theorem 6 be satisfied. Further, suppose that the polynomials
$f, \frac{\partial f}{\partial X_1}, \ldots , \frac{\partial f}{\partial X_i}
\in \mQ [ X_1,\ldots , X_n] $ are given by a straight-line program $\beta$ in 
$\mQ [X_1, \ldots , X_n]$ without essential divisions, and let $ L$ be the 
size of $\beta$. Then there is an arithmetical network with parameters in $\mQ$
 that constructs the following items from the input $\beta$
\begin{itemize}
\item a regular matrix of $\mbox{} \mQ^{n\times n}$ given by its elements that transforms the
variables $X_1, \ldots , X_n$ into new ones $Y_1, \ldots , Y_n$
\item a non-zero linear form $U \in \mQ [Y_{r+1}, \ldots , Y_n]$
\item a division-free straight-line program $\gamma$ in $\mQ[Y_1, \ldots ,Y_r, U]$
that represents non-zero polynomials $\varrho \in \mQ[Y_1, \ldots , Y_r] $ and
$q,p_1,\ldots , p_n \in \mQ [Y_1,\ldots , Y_r, U]$.
\end{itemize}
These items have  the following  properties:
\begin{itemize}
\item[(i)] The variables $Y_1,\ldots ,Y_n$ are in Noether position with respect to
the variety $W^\ast_{n-r}$, the variables $Y_1, \ldots , Y_r$ being free
\item[(ii)] The non-empty open part $(W^\ast_{n-r})_\varrho$ is defined by the
ideal $(q, \varrho X_1-p_1, \ldots ,$
 $\varrho X_n-p_n)_\varrho$ in the localization
$\mQ[X_1,\ldots ,X_n]_\varrho$ .
\item[(iii)] The polynomial $q$ is monic in $u$ and its degree is equal to\linebreak 
$\delta^\ast_{n-r}= \deg^\ast W_{n-r}= \deg W^\ast_{n-r} \le \delta^\ast$.
\item[(iv)] $\max \{ \deg_up_k | 1\le k \le n \} < \delta^\ast_{n-r},\quad
\max \{ \deg p_k | 1\le k \le n \} = (d  \delta^\ast)^{0(1)},$\linebreak
$\deg \varrho = (d \delta^\ast)^{0(1)}$.
\item[(v)] The nonscalar size of the straight-line program $\gamma$ is given 
by $(sd\delta^\ast L)^{0(1)}$.
\end{itemize}
\end{lemma}
The proof of Lemma 9 can be performed in a similar way as in  
\cite{gihemorpar}, \cite{gihemopar} for establishing
the algorithm. For the case handled here, in the i-th 
step one has to apply the algorithm to the localized sequence
$\left( f, \frac{\partial f}{\partial X_1}, \ldots , \frac{\partial f}{\partial X_i} 
\right)_{\Delta}$ as input.
The only point we have to take care of is the process of cleaning 
extranious $\mQ$-irreducible components. Whereas in the proofs of the 
algorithms we refer to it suffices to clean out components lying in a 
prefixed hypersurface (e.g. components at infinity), the cleaning process we
need here is more subtile.
We have to clean out all non-real $\mQ$-irreducible components that appear
during our algorithmic process. The idea of doing this is roughly as follows.
Due to the generic position of the variables $X_1,\ldots,X_n$ all 
$\mQ$-irreducible components of the variety $\widetilde{W}_{n-r}$ can be 
visualized as $\mQ$-irreducible factors of the polynomial $q(X_1,
\ldots,X_r,U)$. If we specialize {\em generically} the variables 
$X_1,\ldots,X_r$ to {\em rational} values $\eta_1,\ldots,\eta_r$, then
by Hilbert's Irreducibility Theorem (in the version of \cite{lang}) 
the $\mQ$-irreducible factors of the {\em multivariate} polynomial
$q(\eta_1,\ldots,\eta_r,U)$ correspond to the $\mQ$-irreducible factors
of the {\em one--variate} polynomial $q(\eta_1,\ldots,\eta_r,U) \in \mQ[U]$.
In order to explain our idea simpler, we assume that we are able to 
choose our specialization of $X_1,\ldots,X_r$ into $\eta_1,\ldots,\eta_r$
in such a way that the hyperplanes $X_1 - \eta_1 = 0,\ldots,X_r - \eta_r = 0$
cut any {\em real} component of $\widetilde{W}_{n-r}$ (This condition is
open in the strong topology and doesn't represent a fundamental restriction
on the correctness of our algorithm. Moreover, our assumption doesn't
affect the complexity). Under these assumptions the $\mQ$-irreducible factors 
of $q(X_1,\ldots,X_r,U)$, which correspond to the real components
of $\widetilde{W}_{n-r}$, reappear as $\mQ$-irreducible factors of 
$q(\eta_1,\ldots,\eta_r,U)$ which contain a real zero. These
$\mQ$-irreducible factors of $q(\eta_1,\ldots,\eta_r,U)$ can be found by a 
factorization procedure and by a real zero test of standard features of 
polynomial complexity character. Multiplying these factors 
and applying to the result the lifting-fibre process of \cite{gihemorpar},
\cite{gihemopar} we find the product $q^*$ of the $\mQ$-irreducible factors of $q(
X_1,\ldots,X_r,U)$, which correspond to the union of the real 
components of the variety $\widetilde{W}_{n-r}$, i.e. to the real part of 
$\widetilde{W}_{n-r}$. The ideal $(q^*,\varrho X_1-p_1,\ldots,\varrho
X_n-p_n)_\varrho$ describes the localization of the real part of 
$\widetilde{W}_{n-r}$ at $\varrho$. All we have pointed out is executible in 
polynomial time if a factorization of univariate polynomials over 
$\mQ$ in polynomial time is available and if our geometric assumptions on
the choice of the specialization is satisfied.

\begin{theorem}
Let the notations and assumptions be as in Theorem 6. Suppose
that the polynomial $f$ is given by a straight-line program
$\beta$ without essential divisions in $\mQ [X_1,\ldots ,X_n]$,
and let $L$ be the nonscalar size of $\beta$.  Further,
let $\delta^\ast_i := \deg^\ast W_i, \; 
\delta^\ast :=  \max \{ \delta^\ast_i | 0
\le i < n \}$ be the corresponding real degrees of the polar
varieties in question, and let $d  :=
\deg f$. Then there is an arithmetical network of size $(n d
\delta^\ast L)^{0(1)}$ with parameters in $\mQ$ which
produces, from the input $\beta$, the coefficients of a non-zero
linear form $u \in \mQ [X_1, \ldots ,X_n]$  and non-zero
polynomials $q,p_1, \ldots , p_n \in \mQ[U]$ showing the following properties:
\begin{enumerate}
\item For any connected component $C$ of $V$ there is a point $\xi \in C$ and 
an element $ \tau \in \R$ such that $q (\tau)=0$ and $\xi = (p_1(\tau), \ldots, p_n(\tau))$
\item $\deg (q) = \delta^\ast_{n-1} \le \delta^\ast$
\item $\max \{ \deg (p_i) | 1 \le i \le n \} < \delta^\ast_{n-1}$.
\end{enumerate}
\end{theorem}

\end{document}